\documentclass[a4paper,10pt]{article}
 \usepackage[margin=1.2in]{geometry}
\usepackage[utf8]{inputenc}
\usepackage{acronym}
\usepackage{color}
\usepackage{amsmath}
\usepackage{amsfonts}
\usepackage{amssymb}
\usepackage{authblk}
\usepackage{mathrsfs}
\usepackage{graphicx}
\usepackage{hyperref}
\usepackage{epstopdf}
\usepackage{graphicx}
 \usepackage[utf8]{inputenc}
 \usepackage[T1]{fontenc}
 \usepackage{lmodern}
 \usepackage[normalem]{ulem}
 \usepackage{verbatim}
 \usepackage{bbm}
 \usepackage{stmaryrd}
 \usepackage{amsmath}
 \usepackage{amssymb}
 \usepackage{dsfont}
\usepackage{amsthm}
\usepackage{hyperref}
\usepackage{relsize}
\usepackage{exscale}

\newtheorem{theo}{Theorem}[section]%

\newtheorem{lemma}[theo]{Lemma}%

\def\tr{{\rm Tr}}

\newcommand{\E}{\mathbb E}

\newcommand{\Pp}{\mathbb P}
\newcommand{\C}{\mathbb C}
\newcommand{\R}{\mathbb R}
\newcommand{\N}{\mathbb N}

\def\la{{\lambda}}
\def\Tr{{\rm Tr}}

\renewcommand\phi\varphi

\title{Large deviations for the largest eigenvalues and eigenvectors of spiked random matrices}

\author[1]{Giulio Biroli}
\author[2]{Alice Guionnet}
\affil[1]{\small Laboratoire de Physique de l'Ecole normale sup\'erieure, ENS, Universit\'e PSL, CNRS, Sorbonne Universit\'e, Universit\'e Paris-Diderot, Sorbonne Paris Cit\'e, Paris, France}
\affil[2]{\small Universit\'e de Lyon, CNRS, ENSL, 46 all\'ee d'Italie, 69007 Lyon}

\begin{document}

\maketitle

\begin{abstract}
We consider matrices formed by a random $N\times N$ matrix drawn from the Gaussian Orthogonal Ensemble (or Gaussian Unitary Ensemble) plus a rank-one perturbation of strength $\theta$, and focus on the largest eigenvalue, $x$, and the component, $u$, of the corresponding eigenvector in the direction associated to the rank-one perturbation. We obtain  the large deviation principle governing the atypical joint fluctuations of $x$ and $u$. Interestingly, for $\theta>1$, in large deviations characterized by a small value of $u$, i.e. $u<1-1/\theta$, the second-largest eigenvalue pops out from the Wigner semi-circle and the associated eigenvector orients in the direction corresponding to the rank-one perturbation.  
We generalize these results to the Wishart Ensemble, and we extend them to the first $n$ eigenvalues and the associated eigenvectors.
\end{abstract}

\section{Introduction}
The large deviations theory for the  spectral properties of random matrix models is a very active domain of research 
in probability theory and theoretical physics.\\
 A lot of works have been devoted to the statistics of the eigenvalues. Following Voiculescu pioneering work on non-commutative entropy \cite{Vo02}, G. Ben Arous and one of the author derived a large deviation principle for the distribution of the empirical measure of the eigenvalues of Gaussian ensembles in the late nineties (in physics known as Coulomb gas method \cite{vivobook}). The proof is based on the explicit density of the joint law of the eigenvalues and the speed of the large deviation principle is the square of the linear dimension of the random matrix.  More than ten years later, C. Bordenave and  P. Caputo \cite{BordCap} obtained a large deviations principle for the same empirical measure but for a Wigner matrix with heavy tails entries, in the sense that their tail decays more slowly than a Gaussian variable at infinity. Their approach is totally different as it follows from the ideas that deviations are created by a few big entries: the rate then depends on the speed of decay of the tail. The large deviation principle for the spectral measure  in the general sub-Gaussian case is still an open problem. \\
 Instead of considering the deviations of the empirical measure, it is also natural to try to understand the probability of deviations of a single eigenvalue.  The deviations of an eigenvalue inside the bulk is closely related to that of the empirical measure but one can seek for the probability of deviations of the extreme eigenvalues. This was achieved for Gaussian ensembles in the Appendix of \cite{BADG01}, see also \cite{Majum}, where it was shown that the large deviations are on the scale of the dimension. Again, the proof was based on the explicit joint law of the eigenvalues. The large deviations principle for the largest eigenvalue was derived in \cite{fanny} for heavy tails. In the case of sharp sub-Gaussian entries, which include Rademacher (binary) entries,  it was recently proved that the large deviations of the extreme eigenvalues are the same than in the Gaussian case \cite{GuHu}. \\
 The probability of atypical eigenvectors has been  much less studied. Again, the only result that we know concerns the Gaussian ensembles: in this case, the invariance by multiplication of the Haar measure implies that each eigenvector is uniformly distributed on the sphere. In \cite{BADG01}, the large deviations for the empirical measure of the properly rescaled entries of an eigenvector was established. The large deviations for the supremum of the entries could also be easily derived.\\ In this article, we address a different question. We want to investigate the large deviations of the eigenvector in a given fixed direction.  In many solvable random matrix models, eigenvectors are uniformly distributed; hence there are no 
meaningful atypical fluctuations or special directions to focus on. 
For a spiked GOE matrix, i.e. a random $N\times N$ matrix drawn from the Gaussian Orthogonal Ensemble plus a rank-one perturbation, there is instead a special direction: the one related to the perturbation. In this case an interesting phenomenon, called BBP-transition, takes place  
by varying the strength of the perturbation (called $\theta$ in the following). As shown in \cite{edwards1976eigenvalue} and then proved rigorously in \cite{baik2005phase} the largest eigenvalue, $x$, pops out of the semi-circle if the perturbation is strong enough. More precisely, $x$ is almost surely equal to two for $\theta\le1$ and to $\theta+1/\theta$ for $\theta>1$. In the latter case, the square of the component of the associated eigenvector 
in the direction associated to the perturbation, that henceforth we shall denote $u$, is almost surely equal to $1-1/\theta^2$. 
In this context the question we raised before becomes meaningful, and it is natural to focus on the good rate function (GRF) that controls the {\it joint} atypical fluctuations of $x$ and $u$.\\
This GRF plays an important role for the geometric properties of random high-dimensional energy landscapes, which can exhibit a number of critical points that is exponentially large in the number of dimensions, as obtained in \cite{cavagna1998stationary,fyodorov2004complexity,bray2007statistics} and rigorously proven and extended in \cite{auffinger2013random,subag2017complexity}. The rigorous method developed to perform those studies is based 
on a large dimensional version of the Kac-Rice formula \cite{adler2009random}, and is strongly related to random 
matrix theory, since the Hessian of the energy function at the critical points---a crucial element in the theoretical analysis---is a random matrix. In order to analyze the dynamics in those rough landscapes 
it is important to know not only the behavior of typical critical points, but also of atypical ones 
associated to index one saddles connecting minima \cite{ros2018complexity}. One has therefore to study 
large deviations of the Hessian, i.e. one needs to condition the critical points
to be of index one and to have the eigenvector associated to the negative eigenvalue oriented in the direction connecting the minima, which leads in fact the problem discussed above.  \\
Noise dressing and cleaning of empirical correlation matrices is another context in which the kind of large deviations addressed in this paper are relevant. In this case, a model that is often considered to interpret the data is the one of spiked Wishart random matrices, whose eigenvalue distribution consists in a Marchenko-Pastur law plus a few eigenvalues that pop out from it. Those few eigenvalues correspond to the signal buried in the noise and the associated eigenvectors play an important 
role in assessing the structure of the correlations, with important applications such as portfolios risk management \cite{bouchaud2000theory}. A natural question in this context is to characterize the joint atypical fluctuations of the largest eigenvalues 
and associated eigenvectors that carry the signal. In this work we obtain the large deviation function that governs them.
\section{Main results}
We consider the matrix
$$Y=X+\theta ww^T$$ where $X$ is from the GOE if $\beta=1$ (resp. the GUE if $\beta=2$) and $\theta$ is a non-negative real number.  $w$ is a fixed unit vector and 
we may assume without loss of generality that $w=e_1=(1,0\cdots,0)$. Let $\lambda_N\ge\lambda_{N-1}\ge\cdots\ge  \lambda_1$ be the eigenvalues of $Y$, with respective eigenvectors $v_N,\ldots, v_1$. The joint large deviations of the largest eigenvalue $\lambda_{N}$ and the component  $|v_{N}(1)|^{2}$ of the associated eigenvector along $w=e_{1}$ is governed by the following theorem.  
\begin{theo}\label{toto1}
The joint law of $ (\lambda_N, |v_N(1)|^2)$ satisfies a large deviation principle in the scale $N$ and good rate function  $I_\beta$. In other words, 
 for any closed set $F$ of $\mathbb R\times [0,1]$ 
\[
\limsup_{N\rightarrow \infty }\frac{1}{N}\log P((\lambda_N, |v_N(1)|^2)\in F)\le -\inf_F I_\beta\,.
\]
and for any open set $O$ of $\mathbb R\times [0,1]$ 
\[
\liminf_{N\rightarrow \infty }\frac{1}{N}\log P((\lambda_N, |v_N(1)|^2)\in F)\ge -\inf_O I_\beta\,.
\]
Moreover, $I_{\beta}$ is a good rate function in the sense that it is non-negative and with compact level sets.
More precisely, the function $I_\beta $ is infinite outside of $S=[2,+\infty)\times [0,1]$ and otherwise given
by
$ I_\beta (x,u)=\beta (I(x,u)-\inf_{S}\{I\})$ where
$$I(x,u)=\frac{x^2}{4}-\frac{1}{2}\theta x u-\int\ln|x-t|d\sigma(t)-\frac{1}{2}\ln |1-u|\qquad\qquad $$
\begin{equation}\label{varpb}
\qquad\qquad -\sup_{2\le y\le x}\{ J(\sigma, \frac{1}{2}\theta(1-u), y)-\frac{y^2}{4}+\int\ln|y-t|d\sigma(t)\}\,.
\end{equation}
where if $\theta\le \frac{1}{2} G_\sigma(x)$, then
$$J(\sigma, \theta,x) = 
\theta^2,$$
whereas if $\theta> \frac{1}{2} G_\sigma(x)$,
$$J(\sigma, \theta,x) = 
\theta x-\frac{1}{2} +\frac{1}{2} \log \frac{1}{2\theta}-\frac{1}{2}\int \log|x-y|d\sigma(y).$$
Here $\sigma(dx)=\sqrt{4-x^2} dx/2\pi$ is the semi-circle distribution and $G_\sigma$ its Cauchy transform.\\
The second largest eigenvalue converges almost surely to the maximizer, $y$, of the variational problem (\ref{varpb}) defined above.  
\end{theo}
We have more explicit results on the rate function, the behavior of the second largest eigenvalue and the component 
of the associated eigenvector along $w=e_{1}$:
\begin{itemize}
\item  {\it For $\theta\ge 1$ and $x> \theta+1/\theta$}: The second eigenvalue pops out of the semicircle for $u<1-1/\theta$, and is equal to $y(u)=\theta(1-u)+\frac{1}{\theta(1-u)}$. The minimum  on $u$
of the large deviation function $I_{\beta}$, for a given $x$, is reached at $u_\theta=1-\frac{\theta x-\sqrt{(\theta x)^2-4\theta^2}}{2\theta^2}$.
\item {\it For $\theta\ge 1$ and $2\le x<\theta+1/\theta$}: The second eigenvalue pops out of the semicircle for $u<1-1/\theta$, and is equal to $\inf(y(u),x)$, i.e. it increases when $u$ decreases until reaching the value $x$. The minimum of the large deviation function is reached at $u_\theta$.
\item {\it For $\theta<1 $ and $x\ge 2$}: The second eigenvalue sticks to two. 
 The minimum of the large deviation function is at $u=0$, if $u_\theta$ is not positive, or at $u_\theta$ otherwise. 
 The latter case corresponds to large enough values of $x$. 
 \item The component $|v_{N-1}(1)|^{2}$ of the eigenvector associated to the second largest eigenvalue is different from 
 zero if and only if the associated eigenvalue is larger than two, i.e. the eigenvector $v_{N-1}$ orients in the direction of 
 $w$ when $\lambda_{N-1}$ pops out from the semi-circle. 
\end{itemize}
An example of the large deviation function (GRF) is shown in Fig. 1  for $\theta=3$.
\\
\begin{figure}[hbtp]
\centering
\includegraphics[width=0.65\columnwidth]{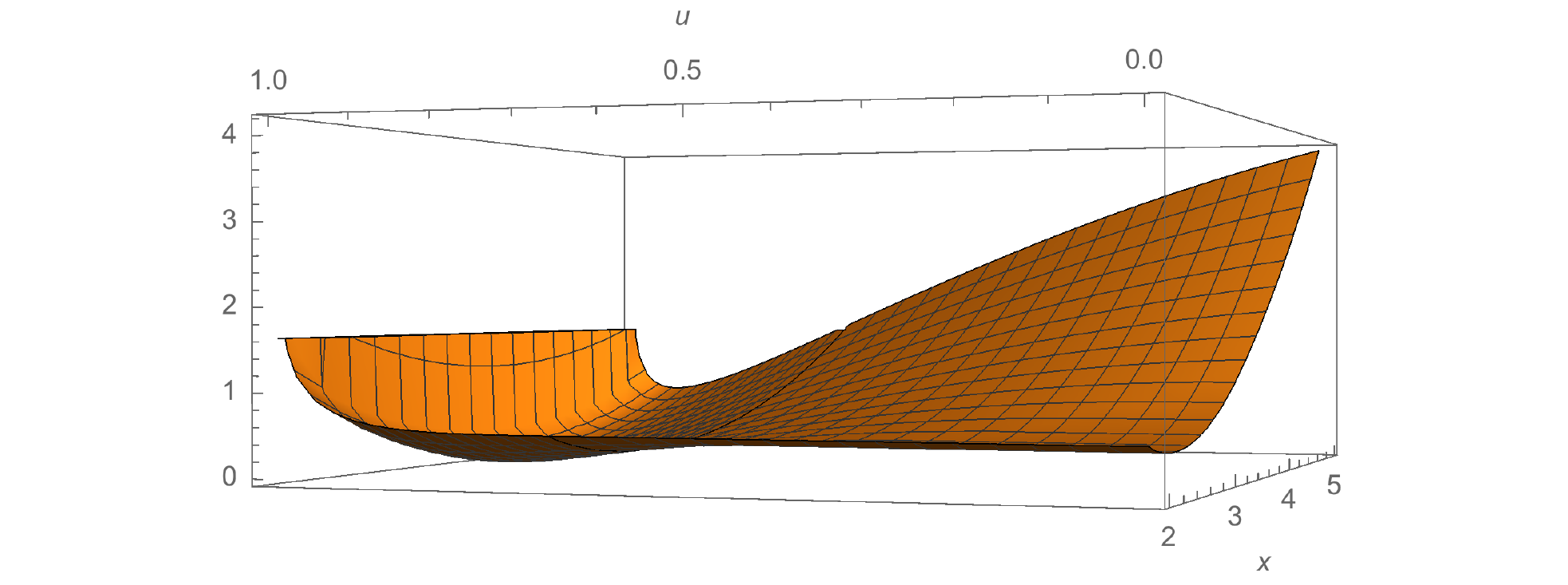}
\caption{Large deviation function plotted for $2<x<5$, $0<u<1$ and $\theta=3$. The global minimum is attained at 
$u=1-1/\theta^2, x=\theta+1/\theta$.}
\label{fig1}
\end{figure}

The previous results can be extended to the large deviations of the $n$ largest eigenvalues $\lambda_{N}\ge\lambda_{N-1}\ge\cdots \lambda_{N-n}$ and their 
first components $|v_{N-i}(1)|^2$. We denote  by $\Delta_{n}=\{x_{i}\in [2,\infty)^{n}:x_{0}\ge x_{1}\ge\cdots \ge x_{n}\}$ the space where these extreme eigenvalues live.
\begin{theo}\label{toto2} Let $n$ be a fixed integer number.
The joint law of $ \{\lambda_{N-i}\}$ and $\{|v_{N-i}(1)|^2\}$, $i=0,...,n-1$, satisfies a LDP with speed $N$ and GRF which is infinite outside of $\Delta_{n}\times [0,1]^n$ and otherwise given
by
$\mathcal I_\beta (\{x_i\},\{u_i\})=\beta (I(\{x_i\},\{u_i\})-\inf\{I(\{x_i\},\{u_i\})\})$ where
$$I(\{x_i\},\{u_i\})=\sum_{i=0}^{n-1}\left(\frac{x_i^2}{4}- \frac{1}{2}\theta x_i u_i-\int\ln|x_i-t|d\sigma(t)\right)$$
$$\qquad \qquad -\sup_{2\le y\le x_{n-1}}\{ J\left(\sigma, \frac{1}{2}\theta\left(1-\sum_{i=0}^{n-1}u_i\right), y\right)-\frac{y^2}{4}+\int\ln|y-t|d\sigma(t)\}-\frac{1}{2}\ln \left|1-\sum_{i=0}^{n-1}u_i\right| $$
where the function $J$ is the same one of Theorem 2.1.
The $n+1$-th largest eigenvalue is equal almost surely to the maximizer $y$ of the variational problem defined above.  
\end{theo}


Moreover, we can as well extend our results in the case of Wishart matrices with covariance which is a finite dimensional perturbation of the identity. To simplify, let us assume it is one dimensional, and consider the Wishart matrix
$$W=\Sigma^{1/2} YY^{*}\Sigma^{1/2}$$
where $Y$ is a $M\times N$ random matrix with i.i.d standard Gaussian entries with variance $1/N$. $\Sigma$ is a $M\times M$ non-negative definite matrix :$\Sigma=I+\gamma ee^{*}$ with $e$ a unit vector.  We assume $M\le N$. We recall that when $M/N$ converges towards $\alpha\in (0,1]$, the empirical measure of $YY^*$ converges towards the so-called Marchenko-Pastur \cite{PM} law $\pi_\alpha$ with support $[1-\sqrt{\alpha}, 1+\sqrt{\alpha}]$. We can study the joint large deviation of the largest eigenvalue $\lambda_N$ and the strength  $|v_N(1)|^2$ of the  eigenvector in the direction $e$ for $W$  as well. We find that 

\begin{theo}\label{toto3}
The joint law of $ (\lambda_N, |v_N(1)|^2)$ satisfies a large deviation principle in the scale $N$ and good rate function  $I^W_\beta$. The function $I_\beta $ is infinite outside of $S=[1+\sqrt{\alpha},+\infty)\times [0,1]$ and otherwise given
by
$ I^W_\beta (x,u)=\frac{\beta}{2} (I^W(x,u)-\inf_{S}\{I^W\})$ where
$$I^W(x,u)=  I(x) -\frac{\gamma}{(1+\gamma)}xu -\alpha\ln(1-u)-\sup_{y\le x}\{2\alpha J(\pi_\alpha,\frac{\gamma(1-u)}{2\alpha (1+\gamma)},y)-I(y)\}$$
where $I(y)= y-(1-\alpha)\ln y-2\alpha\int\ln |y-t|d\pi_\alpha(t)$ is the rate function for the large deviation of the largest eigenvalue of a Gaussian Wishart matrix with covariance equal to the identity.

\end{theo}

\section{Strategy of the proof}
We next focus on the pertubed Wigner matrix $Y=X+\theta ww^{T}$. The law of $Y$ is given by
$$d\mathbb P(Y)=\frac{1}{Z_N} \exp\left\{-\frac{\beta N}{4} \Tr(Y-\theta ww^T)^2\right\} dY=\frac{1}{\tilde Z_N} \exp\left\{-\frac{\beta N}{4} \Tr(Y)^2+\frac{N}{2}\beta\theta  \langle w,Yw\rangle\right\} dY\,.$$
Therefore, since $\langle w,Yw\rangle=\sum\lambda_i | v_i(1)|^2$ when $w=e_{1}$,   the joint law of $(\lambda_N, |v_N(1)|^2)$ is given by
$$dP_N(x,u) =\frac{1}{\bar Z_N} e^{-\frac{\beta N}{4}x^2  +\frac{N\beta}{2}\theta xu}\int \prod_{i=1}^{N-1}|x-\lambda_i|^\beta
d\Pp_{N-1}^{x}(\lambda )\mathbb E[ e^{\frac{N}{2}\beta\theta \sum_{i=1}^{N-1} \lambda_i |v_i(1)|^2}||v_N(1)|^2](u)dx dU_N(u)$$ 
where
$$d\Pp_{N-1}^x(\lambda)=\frac{1}{Z_{N-1}^\infty} \prod_{1\le i<j\le N-1}|\lambda_i-\lambda_j|^\beta  e^{-\frac{\beta N}{4}\sum_{i=1}^{N-1} \lambda_i^2}\prod_{i=1}^{N-1} 1_{\lambda_i\le x}d\lambda_i\,,$$
if
$$Z_{N-1}^\infty=\int \prod_{1\le i<j\le N-1}|\lambda_i-\lambda_j|^\beta  e^{-\frac{\beta N}{4}\sum_{i=1}^{N-1} \lambda_i^2}\prod_{i=1}^{N-1} d\lambda_i\,.$$
$\mathbb E[ .||v_N(1)|^2](u)$ is the expectation on $|v_i(1)|,i\le N-1$ conditionaly to $|v_N(1)|^2=u$.
$dU_N(u)$ is the distribution of $|v_N(1)|^2$.

Our  main goal is to estimate the density of $P_{N}$ when $N$ is large and to apply Laplace's method. 
We infer from concentration inequalities \cite{MM} that with  $\hat\mu^{N-1}=\frac{1}{N-1}\sum_{i=1}^{N-1}\delta_{\lambda_i}$
\begin{equation}\label{conc} \sup_{x>2}\limsup_{N \to \infty} \frac{1}{N} \ln \Pp_{N-1}^x\left( d( \hat \mu^{N-1}, \sigma) > N^{-\kappa'}  \right)\le \limsup_{N \to \infty} \frac{1}{N} \ln \Pp_{N-1}^\infty \left( d( \hat \mu^{N-1}, \sigma) > N^{-\kappa'}  \right)
 = - \infty \,,\end{equation}
We deduce that for $x>2$ (the singularity of the log can be overcome as in \cite{BADG01})
\begin{equation}\label{t1}\prod_{i=1}^{N-1}|x-\lambda_i|=e^{N\int\log|t-x|d\sigma(t)+o(N)}\,,\end{equation}
i.e. we can replace the empirical distribution $\hat\mu^{N-1}$ with its average. 
To estimate the other terms, we first observe that since $(v_i(1))_{1\le i\le N}$ is uniformly distributed on the sphere, we can represent $(v_{i}(1))_{1\le i\le N}$ as
\begin{equation}\label{f1} v_i(1)=\frac{g_i}{(\sum_{i=1}^{N-1} |g_i|^2 +|g_N|^2)^{1/2}}\end{equation}
with independent standard Gaussian variables $g_i,1\le i\le N$, which are real when $\beta=1$ and complex when $\beta=2$.
As a consequence, we see that the distribution $U_{N}$ of $|v_{N}(1)|^{2}$ is
 the Beta-distribution
\begin{equation}\label{t2}dU_N(u)=C_N u^{\beta /2} (1-u)^{(N-1)\beta/2} du\end{equation}

Hence the main point of the proof is to   estimate

\begin{equation}\label{f2}I_N(\lambda,\theta,|v_N(1)|^2)=\mathbb E[ e^{\frac{\beta N}{2}\theta \sum_{i=1}^{N-1} \lambda_i |v_i(1)|^2}||v_N(1)|^2](u)\end{equation}
where the expectation is over $v_i(1),i\le N-1$ with given $|v_N(1)|^2=u$. Thanks to \eqref{f1}, if we fix $v_N(1)$ and denote $g^{N-1}=(g_1,\ldots,g_{N-1})$, we have
$$\|g\|_2^2=\|g^{N-1}\|_2^2 +|v_N(1)|^2 \|g\|_2^2$$
so that  $w_i=g_i/\|g^{N-1}\|_2, 1\le i\le N-1,$ follows the uniform law on the sphere $S^{N-2}$ and 
$$v_i(1)=w_i \frac{\|g^{N-1}\|_2}{\|g\|_2}= \sqrt{1-|v_N(1)|^2} w_i$$
Observe that $(w_i)_{1\le i\le N-1}$ is independent of $|v_N(1)|^2$ (as the computation of the joint law reveals).
Hence
\begin{equation}\label{f3}
I_N(\lambda,\theta,|v_N(1)|^2)=\mathbb E_{w_i, i\le N-1}[ e^{\frac{N}{2}\beta\theta (1-|v_N(1)|^2) \sum_{i=1}^{N-1} \lambda_i |w_i|^2}]\end{equation}

\section{Asymptotic of spherical integrals}
Recall the definition of spherical integrals:
$$I_N(X, \theta)=\E_e[ e^{\beta \theta N\langle e, X e\rangle}]$$
where $e$ is uniformly sampled on the sphere $\mathbb S^{N-1}$ with radius one. 
The asymptotics of 
$$J_N(X,\theta)=\frac{1}{N}\log I_N(X,\theta)$$ 
were studied in \cite{GuMa05} where  the following result was proved.
\begin{theo}\cite[Theorem 6]{GuMa05}\label{myl}
Let $(E_N)_{N \in \N}$ be a sequence of $N \times N$ real symmetric matrices  such that :
\begin{itemize}
\item The sequence of empirical measures $\hat\mu^{N}_{E_N}$ converges weakly to a compactly supported measure $\mu$.
\item There are two real numbers $\lambda_{\min}(E), \lambda_{\max}(E)$ such that 
$$\lim_{N \to + \infty} \lambda_{\min}(E_N) = \lambda_{\min}(E), \quad \lim_{N \to +\infty} \lambda_{\max}(E_N) = \lambda_{\max}(E).$$
\end{itemize}
For any $\theta \geq 0$, 
\[ \lim_{N \to + \infty} J_N(E_N,\theta) =\beta  J(\mu,\theta, \lambda_{\max}(E)) \]
\end{theo}
The limit $J$ is defined as follows.  For a compactly supported probability measure $\mu \in\mathcal{P}(\R)$ we define its Stieltjes transform $G_\mu$ by,
\[ \forall z \in \C \setminus \mathrm{supp}(\mu), \ G_{\mu}(z) := \int_{\R} \frac{1}{z-t} d \mu(t), \]
where $\mathrm{supp}(\mu)$ denotes the support of $\mu$.
In the sequel, for any compactly supported probability measure $\mu,$ we denote by $r(\mu)$ the right edge of the support of $\mu.$ Then  $G_{\mu}$ is a bijection from $(r(\mu), +\infty)$ to $(0, G_{\mu}(r(\mu))$, with
$$G_{\mu}(r(\mu)) = \lim_{t \downarrow r(\mu)} G_{\mu}(t).$$
We denote by $K_{\mu}$ its inverse on $(0, G_{\mu}(r(\mu))$ and let $R_{\mu}(z) := K_{\mu}(z) - 1/z$ be the $R$-transform of $\mu$ as defined by Voiculescu in \cite{Vo5} (defined on $(0, G_{\mu}(r(\mu))$).

In order to define the rate function, we now introduce, for any $\theta \ge 0,$   and $\lambda \ge r(\mu)$,
$$
J( \mu, \theta, \la):= \theta v(\theta, \mu, \la) -\frac{1}{2} \int \log\left(1+2 \theta v(\theta, \mu, \la) - 2\theta y\right) d\mu( y),
$$
with 
$$
 v(\theta, \mu, \la) :=\begin{cases}
                                 R_\mu(2 \theta) & \text{ if } 0 \le 2 \theta\le G_{\mu}( \lambda),\\
                                \lambda - \frac{1 }{2\theta} & \text{ if }2 \theta>  G_{\mu}( \lambda).
                                \end{cases}
$$
In the case where $\mu = \sigma$, the semi-circular law, then,
$$ \forall x> 2, \ G_{\sigma}(x) = \frac{1}{2}(x-\sqrt{x^2-4}), \ R_{\sigma}(x) = x.$$
Therefore,
\begin{lemma}If $\theta\le \frac{1}{2} G_\sigma(x)$, then
$$J(\sigma, \theta,x) = 
{\theta^2},$$
whereas if $\theta> \frac{1}{2} G_\sigma(x)$,
$$J(\sigma, \theta,x) = 
\theta x-\frac{1}{2} +\frac{1}{2} \log \frac{1}{2\theta}-\frac{1}{2}\int \log|x-y|d\sigma(y).$$

\end{lemma}

\section{Proof of Theorem \ref{toto1}}Remark that Theorem \ref{toto1} implies the weak large deviation principle which states that  for $\delta$ small enough,
$$P(|\lambda_N-x|\le\delta+||v_N(1)|^2-u|\le \delta)\simeq e^{-N (I_\beta(x,u)+o(\delta))}\,.$$
Indeed, the weak large deviation principle is simply the restriction of the full large deviation principle to small balls.
To recover the full large deviation principle from its weak version, it is enough to show that the probability is exponentially tight in the sense that deviations mostly occur in a compact set. The latter is easy to check since $|v_N(1)|^2$ lives in a compact set and
$$|\lambda_N|\le \|X\|_\infty +\theta$$
where it is known that $\|X\|_\infty\le M$ with probability greater than $1-e^{-Nc M^2}$ \cite{BADG01}. We refer the reader to \cite{DZ} for more details. Hence, we only need to prove the weak large deviation principle, that is estimate the probability that $(\lambda_{N}, |v_{N}(1)|^{2})$ is close to some $(x,u)$.\\

By Theorem \ref{myl}, if we assume addtionally that $\lambda_{N-1}$ is close to y,  and $|v_N(1)|^2$ close to $u$, we deduce using \eqref{conc} and \eqref{f3} that 
$$I_N(\lambda,\theta,|v_N(1)|^2)\simeq  e^{N\beta  J(\sigma,\frac{1}{2} \theta(1-u), y)+o(N)}$$
But $\lambda_{N-1}$ satisfies a LDP under $\Pp_{N-1}^{x}$ \cite{BADG01} with good rate function which is infinite above $x$ and below $2$, and otherwise given by
$$\beta\left(\frac{y^2}{4}-\int\ln|y-t|d\sigma(t)-C\right)$$
where $C=\inf\{\frac{x^2}{4}-\int\ln|x-t|d\sigma(t)\}$.
Hence we deduce by continuity of the limiting spherical integrals \cite{Ma07} that
$$\int 
d\Pp_{N-1}^{x}(\lambda )\mathbb E[ e^{\frac{\beta N}{2}\theta \sum_{i=1}^{N-1} \lambda_i v_i(1)^2}||v_N(1)|^2]
\simeq e^{N\beta  \sup_{2\le y\le x}\{ J(\sigma, \frac{1}{2}\theta(1-|v_N(1)|^2), y)-\frac{y^2}{4}+\int\ln|y-t|d\sigma(t)+C\}}$$
and therefore, plugging \eqref{t1},\eqref{t2} in the above estimate,
we deduce that the joint law of $(\lambda_N,|v_N(1)|^2)$  is approximately given by
$$
dP_N(x, u)\simeq \frac{1}{Z_N'} \exp\{-\frac{N\beta }{4}x^2  +\frac{N}{2}\beta \theta xu+\beta N\int\ln|x-t|d\sigma(t)\qquad\qquad\qquad$$
$$\qquad \qquad\qquad
+N\beta  \sup_{2\le y\le x}\{ J(\sigma, \frac{1}{2}\theta(1-u), y)-\frac{y^2}{4}+\int\ln|y-t|d\sigma(t)-C\}\} (1-u)^{\beta (N-1)/2}  u^{\frac{\beta}{2}}du dy  $$ 
The final result follows by Laplace's method.

\section{Proof of Theorem \ref{toto2}}
The law of $Y$ is given by
$$d\mathbb P(Y)=\frac{1}{Z_N} \exp\left\{-\frac{\beta N}{4} \Tr(Y-\theta vv^T)^2\right\} dY=\frac{1}{\tilde Z_N} \exp\left\{-\frac{\beta N}{4} \Tr(Y)^2+\frac{\beta N}{2}\theta  \langle v,Yv\rangle\right\} dY$$
and therefore, since $\langle v,Yv\rangle=\sum\lambda_i v_i(1)^2$,   the joint law of $\{\lambda_{N-i},v_{N-i}(1)^2\}$ for $i=0,\dots,n-1$ is given by
$$dP_N^{n}(x,u) =\frac{1}{\bar Z_N} e^{-\frac{\beta N}{4}\sum_{i=0}^{n-1}x_i^2  +\frac{\beta N}{2}\theta \sum_{i=0}^{n-1}x_iu_i}I_{N}^{n}(x,\theta, u) \prod_idx_i dV_N(\{u_i\})$$ 
where $I_{N}^{n}(x,\theta, u)$ equals

$$\int \prod_{i=0}^{n-1}(\prod_{k=1}^{N-n}|x_i-\lambda_k|^\beta) \prod_{i<j}^{n-1}|x_i-x_j|^\beta
d\Pp_{N-n}^{\{x_i\}}(\lambda )\mathbb E[ e^{\frac{N}{2}\beta\theta \sum_{i=1}^{N-n} \lambda_{i=1}^{N-n} |v_i(1)|^2}|v_{N-i}(1)](u_{1},..,u_{n})\,.$$
Here we have denoted
$$d\Pp_{N-n}^{\{x_i\}}(\lambda)=\frac{1}{Z_{N-1}^\infty} \prod_{1\le i<j\le N-n}|\lambda_i-\lambda_j|^\beta  e^{-\frac{\beta N}{4}\sum_{i=1}^{N-n} \lambda_i^2}\prod_{i=1}^{N-n} 1_{\lambda_i\le x_{n-1}}d\lambda_i\,,$$
if
$$Z_{N-n}^\infty=\int \prod_{1\le i<j\le N-n}|\lambda_i-\lambda_j|^\beta  e^{-\frac{\beta N}{4}\sum_{i=1}^{N-n} \lambda_i^2}\prod_{i=1}^{N-n} 1_{\lambda_{i}\le x_{n-1}}d\lambda_i\,.$$
$\mathbb E[ .|\{|v_{N-i}(1)|^2=u_i\}]$ is the expectation on $v_i(1),i\le N-n$ conditionally to $\{v_{N-i}(1)^2=u_i\}$.
$dV_N(\{u_i\})$ is the distribution of $\{v_{N-i}(1)^2\}$. \\
The analysis of the expressions above allows to extend straightforwardly Theorem \ref{toto1} to \ref{toto2}. 
This follows from four remarks:
\begin{enumerate}
\item The term $\prod_{i<j}^{n-1}|x_i-x_j|$ does not lead to any contribution to the GRF as long as 
the $x_i$s are distinct. 
\item The spherical integral $\mathbb E[ e^{\frac{N}{2}\beta\theta \sum_{i=1}^{N-n} \lambda_{i=1}^{N-n} |v_i(1)|^2}|\{v_{N-i}(1)^2=u_i\}]$ is performed on $v_i(1)$ that are uniform on  the sphere $S^{N-n-1}$
and such that $$\sum_{i=1}^{N-n} |v_i(1)|^2+\sum_{i=0}^{n-1}u_i=1$$
As a consequence, one can write 
$$v_i(1)= w_i\sqrt{1-\sum_{i=0}^{n-1}u_i} $$
where $(w_i)_{1\le i\le N-n}$ is uniform on the sphere $S^{N-n-1}$. Therefore the 
spherical integral is the same one evaluated for Theorem \ref{toto1} with $u$ replaced by $\sum_{i=0}^{n-1}u_i$.
\item Because of rotational symmetry $dV_N(\{u_i\})$ depends only on $\sum_{i=0}^{n-1}u_i$. Moreover, calling 
$d\tilde V_N(\sum_{i=0}^{n-1}u_i)$ the distribution of $\sum_{i=0}^{n-1} |v_{N-i}(1)|^2$, we remark that 
$dV_N(\{u_i\})$ and $d\tilde V_N(\sum_{i=0}^{n-1}u_i)$ have the same GRF, since the integral over 
the $u_i$s at fixed $\sum_{i=0}^{n-1}u_i$ does not lead to a term exponential in $N$. Furthermore, 
$d\tilde V_N(\sum_{i=0}^{n-1}u_i)$ and the distribution $dU_N(u)$, introduced for Theorem \ref{toto1}, also 
have the same GRF. In conclusion, the contribution to the total GRF due to $dU_N(\{u_i\})$ is the same 
one of Theorem \ref{toto1} with $u$ replaced by $\sum_{i=0}^{n-1}u_i$.
\item The above implies by Laplace method the weak large deviation principle at any strictly ordered sequence of points $x_0>\cdots>x_n$, that is
$$\limsup_{\delta\rightarrow 0}\limsup_{N\rightarrow \infty}\frac{1}{N}\ln P_N(\cap_{0\le i\le n}\{|\lambda_i-x_i|+|v_N(i)^2-u_i|\le\delta\})=-I(\{x_i\},\{u_i\})\,,$$
and the same when the limsup are replaced by a liminf.
We deduce the same case for $x_i$ ordered and eventually equal by taking approximating sequences $\tilde x_i^\delta$ and $\bar x_i^\delta$ which are strictly ordered and such that 
\begin{eqnarray*}
&&P_N(\cap_{0\le i\le n}\{|\lambda_i-\tilde x^\delta_i|+||v_N(i)|^2-u_i|\le\delta/2\})\\
&&\qquad\le P_N(\cap_{0\le i\le n}\{|\lambda_i-x_i|+||v_N(i)|^2-u_i|\le\delta\})\\
&&\qquad\qquad \le P_N(\cap_{0\le i\le n}\{|\lambda_i-\bar x^\delta_i|+||v_N(i)|^2-u_i|\le 2\delta\})\end{eqnarray*}
By the previous bounds we deduce that
$$-I(\{\tilde x_i^\delta\},\{u_i\})\le \limsup_{N\rightarrow \infty}\frac{1}{N}\ln P_N(\cap_{0\le i\le n}\{|\lambda_i-x_i|+|v_N(i)^2-u_i|\le\delta)\le -\inf_{B^\delta} I(\{y_i\},\{u_i\})$$
where $B^\delta=\cap_{0\le i\le n}\{|y_i-\bar x^\delta_i|+||v_N(i)|^2-u_i|\le 2\delta\}$. The continuity of $I$ in the $x_i$ allows to conclude by letting $\delta$ going to zero.
\end{enumerate}
\section{Study of the rate function}

We can give a more explicit formula of the rate function by noticing that the supremum 
$$ H(x,\theta(1-u))=\sup_{2\le y\le x}\{ J(\sigma, \frac{1}{2}\theta(1-u), y)-\frac{y^2}{4}+\int\ln|y-t|d\sigma(t)\}$$
was already studied in \cite{Ma07}. In the notations of \cite{Ma07}, we are maximizing $-F^1_{\frac{1}{2} \theta(1-u)}(\frac{y}{\sqrt{2}})$ on $y\in [2,x]$. According to \cite[Section 3.2]{Ma07}  of this paper 
we find that
\begin{itemize}
\item If $\theta(1-u)\ge 1$ , the maximum is achieved at $$y(u)=\theta(1-u)+\frac{1}{\theta(1-u)}\,,$$
or at $x$, if $x$ is smaller than $y(u)$, and 
$$H(x,\theta(1-u))= -F^1_{  \frac{1}{2}\theta(1-u)}( \inf(y(u),x)).$$
\item If $\theta(1-u)\le 1$, we are optimizing a decreasing function and therefore the maximum is taken at $2$ and 
$$H(x,\theta(1-u))=J(\sigma, \frac{1}{2}\theta(1-u), 2)+C'=\frac{\theta^2(1-u)^2}{4}+C'$$
with $C'=-1+\int \ln|2-x|d\sigma(x)\,.$
\end{itemize}
Note that these two cases correspond to different asymptotic behaviors of $\lambda_{N-1}$ when $\lambda_N$ goes to $x$ and $|v_N(1)|^2$ goes to $u$: if $\theta(1-u)$ is larger than $1$, $\lambda_{N-1}$ goes to $\inf(y(u),x)$,
and otherwise to $2$.

We can therefore study the optimizer in $u$ of $H$ for a given $x$.
\begin{itemize}
\item For $\theta(1-u)\ge 1$, the contribution to the GRF that depends on $u$ and that we have to minimize reads:
$$-\frac{1}{2}\theta (x-\inf(y(u),x)) u-\frac 1 2 \int\ln|\inf(y(u),x)-t|d\sigma(t)+\frac{\inf(y(u),x) ^2}{4}$$
which is independent of $u$  if $x\le y(u)$. If $x\ge y(u)$ the total derivative of this expression is simply 
$$-\frac{1}{2}\theta (x-\inf(y(u),x))\le 0
$$
since the partial derivative with respect to $y(u)$ is zero ($y(u)$ is an extremum). 
As a consequence, the expression above is a decreasing function of $u$ in the entire available range of $u$.
\item for $\theta(1-u)\le 1$, the contribution to the GRF that depends on $u$ and that we have to minimize reads:
$$-\frac{1}{2}\theta x u-\frac{\theta^2(1-u)^2}{4}-\frac{1}{2}\ln |1-u| $$
We find that the supremum is taken for $\theta>1$ at $u_\theta$ given by
$$1-u_\theta=\frac{\theta x-\sqrt{(\theta x)^2-4\theta^2}}{2\theta^2}$$
which is always smaller than $1/\theta$ for $x\ge 2$, and is non negative  if $x\ge 1+\theta^{-2}$.
For $\theta<1$ the supremum is in $u_\theta$ if it is positive or in $u=0$ otherwise. 
\end{itemize}
\vskip 1cm

Collecting all previous results, we can finally obtain the complete behavior as a function of $u$ for a fixed $x\ge2$:
\begin{enumerate}
\item $\theta\ge 1$ and $x>\theta+1/\theta$. In this case $y(u=0)<x$, when u increases the second eigenvalue $y(u)$ decreases and reaches two for $\theta(1-u)=1$. The minimum takes place for 
$\theta(1-u)<1$ at $u_\theta$.
\item $\theta\ge 1$ and $2<x<\theta+1/\theta$. In this case $y(u=0)>x$, so the LD function remain constant 
until the value $u^*$ for which $y(u)=x$ and then starts decreasing. In this case the second eigenvalue is 
stuck at $x$, it starts decreasing after $u^*$and reaches  two for $\theta(1-u)=1$. The minimum takes place for 
$\theta(1-u)<1$ at $u_\theta$.
\item   $\theta< 1$ and $x>2$. In this case there is only the regime $\theta(1-u)<1$. Depending on the value of 
$x$ the minimum can be at $u=0$ (when $u_\theta$ is not positive) or in $u_\theta$ for $x$ large enough. 
The second eigenvalue remains at two for all $u$. 
\end{enumerate}
Finally, we consider the rate function obtained in Theorem \ref{toto2} for $x_1,x_2$ and the components $u_1$ and $u_2$ of the associated eigenvectors, and study its minimum over $x_2$ and $u_2$. We focus first on the case where the minimum is attained in $x_2>2$ for which 
we know from  Theorem \ref{toto1} that $\theta(1-u_1)>1$.
 Using the results above to obtain explicitly the rate function of $x_1,x_2,u_1,u_2$, and taking the derivative with respect to $u_2$, one finds, for values of $u_2$ such that $\theta(1-u_1-u_2)\ge 1$:
\[
\frac{\partial I}{\partial u_2}=-\frac \theta 2 (x_2-\inf(y(u_1+u_2),x_2))  
\]
This implies that the rate function is flat for $u_2$ such that $y(u_1+u_2)>x_2$ (if any)
and decreases for larger $u_2$, thus implying that $u_2>0$ whenever $x_2>2$. Conversely, if $x_2=2$ then 
the minimum is reached in $u_2=0$. This can be shown by noticing that the contribution to the GRF that depends on $u_2$ is in this case:
$$-\frac{1}{2}\theta x_2 u_2-\frac{\theta^2(1-u_1-u_2)^2}{4}-\frac{1}{2}\ln |1-u_1-u_2| $$
where we have used that $\theta(1-u_1-u_2)<1$ since $\theta(1-u_1)<1$ for $x_2=2$. 
The minimum of this function for $0\le u_2\le 1$ is attained at $u_2=0$.

\section{The case of Wishart matrices}\label{wis}
In the case of $M\times M$ Wishart matrix $W=\Sigma^{1/2} YY^{*}\Sigma^{1/2}$ , we can diagonalize  the matrix $W=UD(\lambda)U^{*}$ where $D(\lambda)$ is a diagonal matrix with entries given by the eigenvalues and $U$ the matrix of the eigenvectors. 
Starting from the law of $L=YY^{*}$ which reads \cite{james}:
$$d\mathbb P_{M}^{L}=\frac{1}{Z_{M}^{L}}e^{-\frac \beta 2 N \tr(L)}\det(L)^{\frac \beta 2 (N-M+1)-1}dL
$$
we obtain the joint law of the eigenvalues and $U$:

$$d\mathbb P_{M}^{W}(\lambda,U)=\frac{1}{Z_{M}^{W}}\prod_{i<j}|\lambda_{i}-\lambda_{j}|^{\beta} e^{-\frac{\beta N}{2}\tr(UD(\lambda)U^{*}\Sigma^{-1})}\prod_{{1\le i\le N}} \lambda_{i}^{\frac{\beta}{2}(N-M+1)-1} 1_{\lambda_{i}\ge 0}d\lambda_{i}dU$$
Writing that $\Sigma^{-1}=I-\frac{\gamma}{\gamma+1} ee^{*}$ we deduce that
$$d\mathbb P_{M}^{W}(\lambda,U)=\frac{1}{Z_{N}^{W}}\prod_{i<j}|\lambda_{i}-\lambda_{j}|^{\beta} e^{-\frac{N\beta }{2}\sum \lambda_{i}+\frac{N\beta \gamma}{2(1+\gamma)} \langle e, UDU^{* }e\rangle}
\prod_{{1\le i\le N}} \lambda_{i}^{\frac{\beta}{2}(N-M+1)-1} 1_{\lambda_{i}\ge 0}d\lambda_{i}dU$$
and we conclude that the joint law of the maximum  eigenvalue $\lambda_M$ and $|\langle e, u_M\rangle|^2$ is given by
$$dP_M^W(x,u)=\frac{1}{\bar Z_M^W} e^{-\frac{\beta N}{2}x+\frac{N\beta \gamma}{2(1+\gamma)}xu} x^{\frac{\beta}{2}(N-M+1)-1} 1_{x\ge 0}dx dU_M(u) I_M(\lambda,\frac{\gamma}{1+\gamma},u)$$
where 
$$I_M(\lambda,\theta,u)
 \int \prod_{i\le M-1}|\lambda_{i}-x|^{\beta}  \mathbb E[ e^{\frac{\beta \theta N}{2} \sum_{i=1}^{M-1} \lambda_i |v_N(i)|^2}||v_N(1)|^2=u] d\mathbb P_{M-1}^x(\lambda)$$
if $\mathbb P_{M-1}^x$ is the law of the remaining $N-1$ eigenvalues of Wishart matrices conditioned to be smaller than $x$. The proof of Theorem \ref{toto3} then follows exactly the same steps as before.\\\\
{\bf Acknowledgements}  We thank Mylene Maida for pointing out to us that the case of Wishart matrices could be handled by our techniques.
GB is partially supported by the Simons Foundation collaboration
Cracking the Glass Problem (No. 454935 to G. Biroli). AG is partially supported by 
 LABEX MILYON (ANR-10-LABX-0070) of Universit\'e de Lyon.

\bibliographystyle{unsrt}

\end{document}